\input amstex
\documentstyle{amsppt}
\magnification=\magstep1 \hsize=5.75 true in \hoffset=.565 true in
\vsize=8.65 true in \voffset=.125 true in \topmatter
\title
THE CRITICAL ORDER OF CERTAIN HECKE $L$-FUNCTIONS OF IMAGINARY
QUADRATIC FIELDS
\endtitle
\rightheadtext{L-functions of imaginary quadratic fields}
\author
Chunlei Liu \& Lanju Xu
\endauthor
\address
Morningside Center of Mathematics, Chinese Academy of Sciences,
Beijing 100080, People's Republic of China
\endaddress
\email chunleiliu\@mail.china.com \& xulanju1012\@sina.com
\endemail
\abstract Let $-D < -4$ denote a fundamental discriminant which is
either odd or divisible by 8, so that the canonical Hecke
character of $\Bbb Q(\sqrt{-D})$ exists. Let $d$ be a fundamental
discriminant prime to $D$. Let $2k-1$ be an odd natural integer
prime to the class number of $\Bbb Q(\sqrt{-D})$. Let $\chi$ be
the twist of the $(2k-1)$th power of a canonical Hecke character
of $\Bbb Q(\sqrt{-D})$ by the Kronecker's symbol
$n\mapsto(\frac{d}{n})$. It is proved that the order of the Hecke
$L$-function $L(s,\chi)$ at its central point $s=k$ is determined
by its root number when $|d| \leq
c(\varepsilon)D^{\frac{1}{24}-\varepsilon}$ or, when $|d| \leq
c(\varepsilon)D^{\frac1{12} -\varepsilon}$ and $k\geq 2$, where
$\varepsilon > 0$ and $c(\varepsilon)$ is a constant depending
only on $\varepsilon$.
\endabstract
\thanks
The first author is supported by NSFC.
\endthanks
\endtopmatter
\document
\head 0. Introduction
\endhead
Let $K$ be an imaginary quadratic field with discriminant $-D<-4$
and class number $h$. Suppose that $D$ is either odd or divisible
by 8. Then, according to Rohrlich [Rc], there are exactly
$\gcd(2,D)h$ Hecke characters $\chi_{can}$ of $K$
satisfying\roster
\item
The conductor of $\chi_{can}$ is $(2\sqrt{-D},D)$, where
$(2\sqrt{-D},D)$ denotes the ideal generated by $2\sqrt{-D}$  and
$D$.
\item
$\chi_{can}((\alpha)) =\pm \alpha$ if $(\alpha)$ is a principle
ideal prime to $(2\sqrt{-D},D)$.
\item
$\chi_{can}((n)) = (\frac{-D}{n})n$  if $n$ is a positive integer
prime to $(2\sqrt{-D},D)$.
\endroster
We call such a character $\chi_{can}$ canonical.Let $d$ be a
fundamental discriminant prime to $D$. Let $2k-1$ be an odd
positive integer prime to the class number $h$. Let $\chi$ be the
product of the $(2k-1)$th power of a canonical Hecke character of
$K$ and the lifting of the Kronecker's symbol
$n\mapsto(\frac{d}{n})$. Let $L(s,\chi)$ be the Hecke $L$-function
attached to $\chi$. Then
$$
\Lambda(s,\chi) = (D^{\ast}|d|
)^{s}(2\pi)^{-s}\Gamma(s)L(s,\chi)=W(\chi)\Lambda(2k-s,\chi),
$$
where $D^{\ast} = D\gcd(2,D)$ and $W(\chi) = \pm1$ is the root
number. It is well known that the Hecke $L$-function $L(s,\chi)$
is the $L$-function of a newform $f$ of level $(D^{\ast}|d|)^2$
and weight $2k$ with coefficients in $\Bbb Q$. Let $M$ be the
Grothendieck motive over $\Bbb Q$ attached to $f$ by U. Jannsen
and A. J. Scholl. According to a conjecture of Beilinson and Bloch
and a result of S. Zhang, the order of the Hecke $L$-function
$L(s,\chi)$ at its central point $s=k$ is closely related to the
arithmetic of $M$. Because of their arithmetic nature, these
$L$-functions have been extensively studied by, among other
authors, Gross [Gr], Rohrlich, Rodriguez-Villegas and Yang.
Throughout this paper $\varepsilon > 0$ is arbitrary small and the
constants $c(\varepsilon)$, $c_1(\varepsilon)$ and
$c_2(\varepsilon)$, and those implied in the symbols $\ll$, $\gg$
and $O$ depend at most on $\varepsilon$ and $k$. We list the
following update results. \roster
\item
If $k = 1$ and $W(\chi) = 1$, Rohrlich [Rb] proved that
$L(1,\chi)$ doesn't vanish when $|d| \ll D^{1/39-\varepsilon}$.
\item
If $k = 1$ and $W(\chi)= -1$, Miller-Yang [MY] proved that
$L'(1,\chi)$ doesn't vanish when  $|d| \ll D^{1/35 -
\varepsilon}$.
\item
If $W(\chi) = 1$, Yang [Ya] proved that $L(k, \chi)$ doesn't
vanish when $D\equiv 7\ (mod\ 8)$, $d\equiv 1\ (mod\ 4)$ is
positive and $\sqrt{D} > d^2(\frac{12}{\pi}\ln d+M(k))$, where
$M(k)$ is a constant depending only on $k$.
\endroster
In this paper, we shall prove the following two theorems.
  \proclaim{Theorem 1}If $W(\chi) = 1$, then
  $L(k,\chi)$ does not vanish when $|d| \ll D^{\frac1{12} -\varepsilon}$.
 \endproclaim
  \proclaim{Theorem 2}If $W(\chi) = -1$, then $L'(k,\chi)$ doesn't vanish when
  $|d| \ll D^{\frac1{24} -\varepsilon}$ or, when $|d| \ll D^{\frac1{12} -\varepsilon}$ and $k\geq 2$.
  \endproclaim
From the theorems one sees that \proclaim{Corollary 1}The order of
the Hecke $L$-function $L(s,\chi)$ at its central point $s=k$ is
determined by its root number if $|d| \ll D^{1/24-\varepsilon}$
or, if $|d| \ll D^{\frac1{12} -\varepsilon}$ and $k\geq 2$.
 \endproclaim
Assume that $k=1$. Then $\chi$ lifts to a Hecke character $\psi$
of the Hilbert class field $H$ of $K$ commuting with the action of
the Galois group of $H$ over $\Bbb Q$. And
$$
L(s,\psi) = \prod_{\phi}L(s,\chi\phi),
$$
where $\phi$ runs over all characters of the class group of $K$.
Let $j$ be the $j$-invariant of an elliptic curve over $H$ with
complex conjugation by the ring of integers of $K$. According to
Gross [Gr], there is a unique elliptic curve $A$ over $H$ with
$j$-invariant $j$ which is isogenous to all its Galois conjugates
and whose $L$-function is
$$L(s,A/H) = L(s,\psi)^2.$$
$A$ descends to two isogenous elliptic curves over $\Bbb Q(j)$
with $L$-function $L(s,\psi)$. Let $A_{d}$ be one of them. By
results of Kolyvagin-Logachev [KL] and Gross-Zagier [GZ], the
theorems imply the following arithmetic consequences.
\proclaim{Corollary 2 }If $W(\chi)=1$ and $|d| \ll
D^{\frac1{12}-\varepsilon}$, then the Mordell-Weil group and the
Shafarevich-Tate group of $A_{d}/\Bbb Q(j)$ are finite.
\endproclaim
\proclaim{Corollary 3}If $W(\chi) = -1$ and $|d| \ll
D^{\frac1{24}-\varepsilon}$, then $A_{d}/\Bbb Q(j)$ has a finite
Shafarevich-Tate group and a Mordell-Weil group of rank $h$.
\endproclaim
{\it Acknowledgements.}This work is motivated by Tonghai Yang's
lectures at Morningside Center of Mathematics, Chinese Academy of
Sciences. The authors thank him for his lectures and for his
suggestions on the manuscript, and thank Shouwu Zhang, Fei Xu and
Kezheng Li for inviting them to visit the center. The first author
wishes to express his thanks to Chengbiao Pan for discussions on
Dirichlet $L$-functions, and to Yuan Wang for drawing his
attention to character sums in algebraic number fields. \head
 1. $L$-functions with root number $1$
 \endhead
In this section we shall prove Theorem 1. Put
$$
L(s,\chi,p) = \sum_{(\alpha)\in p}\chi((\alpha))(N(\alpha))^{-s},
\quad \Re s
> 3/2,
$$
where $p$ is the set of all principal integral ideals of $K$.
According to Rohrlich [Rb], for all ideal class characters
$\varphi$ of $K$, the $L$-functions $L(s,\chi \varphi)$ satisfy
the same functional equation as $L(s,\chi)$ does. So
$$
\Lambda(s,\chi,p) =
(D^{\ast}|d|)^{s}(2\pi)^{-s}\Gamma(s)L(s,\chi,p)=
W(\chi)\Lambda(2k-s,\chi,p).
$$
Note that
$$
L(s,\chi,p) = L_D(2s-2k+1) + L(s),
$$
where
$$
L_D(s)=\sum_{n\geq 1,(n,d)=1}(\frac{-D}{n})n^{-s},
$$
and
$$
L(s) = \sum_{(\alpha)\in p'} \chi((\alpha))  N(\alpha)^{-s}
$$
with $p'$ denoting the set of all principal integral ideals not
generated by rational numbers. It is easy to see that $L_D(s)$ is
the Dirichlet $L$-function attached to the Dirichlet character:
$$
({\Bbb Z}/(D|d|))^{\times}\rightarrow {\Bbb C}^{\times},\quad n\mapsto (\frac{-D}{n}).
$$
We now suppose that $W(\chi) = 1$ and proceed to prove Theorem 1.
According to Shimura [Sh] and Rohrlich [Rb], $L(k,\chi) = 0$
implies $L(k,\chi \varphi^{2k-1}) = 0$ for all ideal class
characters $\varphi$ of $K$. So it also implies $L(k,\chi,p)=0$ as
$$
h L(s,\chi,p) = \sum_{\varphi}L(s,\chi \varphi) =
\sum_{\varphi}L(s,\chi \varphi^{2k-1}),
$$where $\varphi$ runs over all ideal class characters of $K$.Hence, to prove $L(k,\chi) \neq 0$,
it suffices to prove that $L(k,\chi,p) \neq 0$. It follows from
the functional equation and a formula of Cauchy that
$$
\frac{1}{2}L(k,\chi,p) = \frac1{2\pi i}\int_{(2k)}\frac{(2\pi)^{k}
\Lambda(s,\chi,p)}{(D^{\ast}|d|)^{k}\Gamma(k)}\frac{ds}{s-k}.
$$
As
$$
L(s,\chi,p) =L_D(2s-2k+1) + L(s).
$$
we get the following approximation to the central value:
$$
\frac{1}{2}L(k,\chi,p) =I_1 + I_2,
$$
where
$$
I_1 = \frac1{2\pi i} \int_{(2k)} (D^{\ast}|d|)^{s-k} (2\pi)^{k-s}
\frac{\Gamma(s)}{\Gamma(k)} L_D(2s+1-2k) \frac{ds}{s-k}
$$
and
$$
I_2 = \frac1{2\pi i} \int_{(2k)} (D^{\ast}|d|)^{s-k} (2\pi)^{k-s}
\frac{\Gamma(s)}{\Gamma(k)} L(s) \frac{ds}{s-k}.
$$
Theorem 1 now follows from the estimate
 $$
I_2 \ll D^{-\frac1{16}+\varepsilon}|d|^{\frac34+\varepsilon},
 $$
which will be proved in the next section, and the estimate
$$
I_1 \geq c_1(\varepsilon)(D|d|)^{-\varepsilon}-c_2(\varepsilon)(D|d|)^{-\frac1{16}+\varepsilon},
$$
which we are going to prove. Shifting the line of integration in
$I_1$ to $\Re s = k - 1/4$, we get
$$
I_1 = L_D(1) +\frac1{2\pi i} \int_{(k-1/4)} (D^{\ast}|d|)^{s-k}
(2\pi)^{k-s} \frac{\Gamma(s)}{\Gamma(k)} L_D(2s+1-2k)
\frac{ds}{s-k}.
$$
Applying Burgess' estimate [Bu]
$$
L_D(\frac12+ it) \ll (D|d|)^{\frac3{16}+\varepsilon }(|t|+1),
$$
we get
$$
I_1 = L_D(1) + O((D|d|)^{-\frac1{16}+\varepsilon}).
$$
The required estimate for $I_1$ now follows from Siegel's estimate
$$
L_D(1) \gg (D|d|)^{-\varepsilon}.
$$
\head 2. The complex part of the approximation to the central value
\endhead
In this section, we shall prove that
 $$
 I_2 \ll D^{-\frac1{16}+\varepsilon}|d|^{\frac34+\varepsilon}.
 $$
Recall that
$$
I_2 = \frac1{2\pi i} \int_{(2k)} (D^{\ast}|d|)^{s-k} (2\pi)^{k-s}
\frac{\Gamma(s)}{\Gamma(k)} L(s) \frac{ds}{s-k}.
$$
Applying Mellin's inversion and writing
$$
\chi((\alpha)) = \epsilon(\alpha)\alpha^{2k-1},
$$
where $\epsilon$ is a quadratic character with conductor
$d(2\sqrt{-D},D)$ on the subgroup of $K^\times$ consisting of
elements prime to $d(2\sqrt{-D},D)$, we get
$$
 I_2 =
 \frac1{\Gamma(k)}\sum_{(\alpha)\in p'} \epsilon(\alpha)\alpha^{2k-1}N(\alpha)^{-k}
 \int_{\frac{2\pi N(\alpha)}{D^{\ast}|d|}}^{\infty}
 e^{-\xi} \xi^k \frac{d\xi}{\xi}.
 $$
The contribution from the terms with $\Re \alpha$ or $|\Im \alpha
|\geq (D|d|)^{\frac12}\log (D|d|)$ is bounded by
$$ \frac{2}{\Gamma(k)}\sum_{n\geq D^\ast|d|\log (D|d|)} \int_{\frac{2\pi n}{D^{\ast}|d|}}^{\infty}
 e^{-\xi} \xi^k \frac{d\xi}{\xi}
 \ll (D^\ast|d|)^{-1}.
 $$
The subsum from the terms of $I_2$ with $0 < \Re \alpha, |\Im
\alpha |< (D|d|)^{\frac12}\log (D|d|)$ equals
$$
\sum_{u,v}\frac {2(u +
\sqrt{-D}v)^{2k-1}}{\Gamma(k)(u^2+Dv^2)^k}\epsilon(\frac{u +
\sqrt{-D}v}2)
 \int_{\frac{\pi (u^2+Dv^2)}{2D^{\ast}|d|}}^{\infty}
 e^{-\xi} \xi^k \frac{d\xi}{\xi},
 $$
 where $(u,v)$ runs over pairs of integers satisfying
 $$
 0< u, \sqrt{D}|v|< 2(D|d|)^{\frac12}\log (D|d|),\quad
 4|(u^2 + Dv^2).
 $$
Conjugate terms grouped together, it becomes
$$
\sum_{u,v}a(u,v)\epsilon(\frac{u + \sqrt{-D}v}2)\int_{\frac{\pi (u^2+Dv^2)}{2D^{\ast}|d|}}^{\infty}
 e^{-\xi} \xi^k \frac{d\xi}{\xi},
$$
where $(u,v)$ runs over pairs of integers satisfying
 $$
 0< u, \sqrt{D}v< 2(D|d|)^{\frac12}\log (D|d|),\quad
 4|(u^2 + Dv^2),
 $$
and
$$
\frac12 a(u,v) = \frac{(u + \sqrt{-D}v)^{2k-1}+(u -
\sqrt{-D}v)^{2k-1}}{\Gamma(k)(u^2+Dv^2)^k}.
 $$
It splits dyadically into at most $4\log_2^2(2D|d|)$ sums of the
form
$$
\sum_{N\leq v<N'}\sum_{M\leq u<M',4|(u^2 +
Dv^2)}a(u,v)\epsilon(\frac{u + \sqrt{-D}v}2)\int_{\frac{\pi
(u^2+Dv^2)}{2D^{\ast}|d|}}^{\infty}
 e^{-\xi} \xi^k \frac{d\xi}{\xi},
$$
where $0< M, \sqrt{D}N< (D|d|)^{\frac12}\log (D|d|),$ $N'\leq 2N,$
and $M'\leq 2M$. By Abel's summation formula, the inner sum is
bounded by
$$
O(\log^{2k}(D|d|))\min (M^{-1},D^{-1/2}N^{-1})\max_{M<w\leq
2M}|S_v(w)|,
$$
where
$$
S_v(w) = \sum_{M\leq u<w, 4|(u^2 + Dv^2)}\epsilon(\frac{u +
\sqrt{-D}v}2).
$$
We claim that
$$
S_v(w) \ll
|d|+M^{\frac12}D^{\frac3{16}+\varepsilon}|d|^{\frac12},\quad
w<2M,$$ from which the estimate for $I_2$, which is stated at the
beginning of this section, follows.
 Write $\epsilon =\epsilon_0\epsilon_1,$
 where $\epsilon_0$ and $\epsilon_1$ have conductors $\frac{\sqrt{-D}}{(\sqrt{-D},4)}$
and $d(\sqrt{-D},4)$ respectively. Let $k_0$ and $k_1/2$ be the
least positive integers in $\frac{\sqrt{-D}}{(\sqrt{-D},4)}$
and$d(\sqrt{-D},4)$ respectively. Then
$$
S_v(w)= \sum_{1\leq j<k_1, 4|(j^2 + Dv^2)}\sum_{M\leq u<w, u\equiv
k_0j(k_1)} \epsilon_0(u/2)\epsilon_1(\frac{u+ \sqrt{-D}v}2).
$$
The inner sum equals
$$
\epsilon_0(k_1/2)\epsilon_1(\frac{k_0j + \sqrt{-D}v}2)\sum_{\frac{M-k_0j}{k_1} \leq l < \frac{w-k_0j}{k_1}}\epsilon_0(l),
$$which, according to Burgess [Bu], is bounded by $O(1+(M/|d|)^{1/2}D^{\frac3{16} + \varepsilon})$.
So
$$
S_v(w) \ll |d|+M^{\frac12}D^{\frac3{16}+\varepsilon}|d|^{\frac12}
$$
as claimed. \head 3. Twists of root number $-1$
\endhead
In this section we shall prove Theorem 2. Similarly, it suffices
to show that
  $L'(k,\chi,p) \neq 0$ under the condition of Theorem 2 (see [GZ] or [MY]). Suppose that $W(\chi) = -1$.
 It follows from this functional equation and a formula of Cauchy that
$$
\frac12 L'(k,\chi,p) =\frac1{2\pi i} \int_{(2k)}\frac{(2\pi)^{k}
\Lambda(s,\chi,p)}{(D^{\ast}|d|)^{k}\Gamma(k)} \frac{ds}{(s-k)^2}.
$$
As
$$
L(s,\chi,p) = L_D(2s-2k+1) + L(s),
$$
 we get the following approximation to the central derivative:
$$
\frac12L'(k,\chi,p) = R_k + C,
$$
where
$$
R_k = \frac{1}{2\pi i}
\int_{(2k)}(\frac{D^\ast|d|}{2\pi})^{s-k}\frac{\Gamma(s)}{\Gamma(k)}L_D(2s-2k+1)\frac{ds}{(s-k)^2}
$$
and
$$
C = \frac1{2\pi i} \int_{(2k)} (D^{\ast}|d|)^{s-k} (2\pi)^{k-s}
\frac{\Gamma(s)}{\Gamma(k)} L(s) \frac{ds}{(s-k)^2}.
$$
Theorem 2 now follows from the estimate
$$
R _k\geq .0351-c(\varepsilon) (D|d|)^{-\frac1{16}+\varepsilon},
$$
which will be proved in the next section, and the estimate
$$
C\ll \cases D^{-\frac1{16}+\varepsilon}|d|^{\frac32+\varepsilon}, &k = 1,\\
D^{-\frac1{16}+\varepsilon}|d|^{\frac34+\varepsilon}, &k >
1,\endcases$$ which we are going to prove. Applying Mellin's
inversion we get
$$C = \sum_{(\alpha)\in p'}\chi((\alpha))N(\alpha)^{-k}
 \int_{\frac{2\pi N(\alpha)}{D^{\ast}|d|}}
 ^{\infty}e^{-\xi} \xi^k(\log\xi - \log\frac{2\pi N(\alpha)}{D^{\ast}|d|})\frac{d\xi}{\xi}.
$$
The contribution from the terms with $ \Re \alpha$ or $ |\Im
\alpha| \geq (D^\ast|d|)^\frac12\log (D|d|)$ is bounded by $$
\sum_{n\geq D^\ast|d|\log (D|d|)}\int_{\frac{2\pi n}{D^{\ast}|d|}}
 ^{\infty}e^{-\xi} \xi^k(\log\xi - \log\frac{2\pi n}{D^{\ast}|d|})\frac{d\xi}{\xi} \ll(D^\ast|d|)^{-1}.
 $$
The subsum from the terms of $C$ with $0 <\Re \alpha, |\Im \alpha| < (D|d|)^\frac12\log (D|d|)$ equals
$$
 \sum_{u,v}\frac {2(u + \sqrt{-D}v)^{2k-1}}{(u^2+Dv^2)^k}f(u,v)\epsilon(\frac {u + \sqrt{-D}v}2),
 $$
 where $(u,v)$ runs over pairs of integers satisfying
 $$
 0< u, \sqrt{D}|v| < 2(D|d|)^{\frac12}\log (D|d|),\quad
 4|(u^2 + Dv^2)
 $$
 and
 $$f(u,v) =
\int_{\frac{\pi (u^2+Dv^2)}{2D^{\ast}|d|}}
 ^{\infty}e^{-\xi} \xi^k(\log\xi - \log\frac{\pi (u^2+Dv^2)}{2D^{\ast}|d|})\frac{d\xi}{\xi}.
 $$
Conjugate terms grouped together, it becomes
$$
 \sum_{u,v}a(u,v)f(u,v)\epsilon(\frac {u + \sqrt{-D}v}2),
 $$
 where $(u,v)$ runs over pairs of integers satisfying
 $$
 0< u, \sqrt{D}v < 2(D|d|)^{\frac12}\log (D|d|),\quad
 4|(u^2 + Dv^2).
 $$
 It splits dyadically into at most $4\log_2^2(2D|d|)$ sums of
the form
$$
\sum_{N\leq v<N'}\sum_{M\leq u<M',4|(u^2 + Dv^2)}a(u,v)f(u,v)\epsilon(\frac{u + \sqrt{-D}v}2),
$$
where $0< M, \sqrt{D}N< (D|d|)^{\frac12}\log (D|d|),$ $N'\leq 2N,$
and $M'\leq 2M$. By Abel's summation formula, the inner sum is
bounded by
$$
O(D|d|\log^2(D|d|))\min(M^{-3},D^{-3/2}N^{-3})\max_{M<w\leq 2M}|S_v(w)|,
$$
if $k = 1$, and by
 $$
O(\log^{2k}(D|d|))\min(M^{-1},D^{-1/2}N^{-1})\max_{M<w\leq 2M}|S_v(w)|,
$$
if $k > 1$, where
$$
S_v(w) = \sum_{M\leq u<w, 4|(u^2 + Dv^2)}\epsilon(\frac{u
+\sqrt{-D}v}2).
$$
In \S2, we have proved that
$$
S_v(w) \ll
|d|+M^{\frac12}D^{\frac3{16}+\varepsilon}|d|^{\frac12},\quad w<2M.
$$
The desired estimate for $C$ now follows. \head 4. The rational
part of the approximation to the central derivative
\endhead
In this section we shall prove that
$$
R_k \geq .0351-c(\varepsilon) (D|d|)^{-\frac1{16}+\varepsilon}.
$$
Recall that
$$
R_k = \frac{1}{2\pi i}
\int_{(2k)}(\frac{D^\ast|d|}{2\pi})^{s-k}\frac{\Gamma(s)}{\Gamma(k)}L_D(2s-2k+1)\frac{ds}{(s-k)^2}.
$$
A change of variable yields
$$
R_k =\frac{1}{\pi i}
\int_{(k+1)}(\frac{D^\ast|d|}{2\pi})^{s-1}\frac{\Gamma(s+k-1)}{\Gamma(k)}\frac{L_D(2s-1)ds}{(s-1)^2}.
$$
Shifting the line of integration to $\Re s = 3/4$ and applying
Burgess' estimate, we get
$$
R_k = \Lambda'_k(1)+O((D|d|)^{-\frac{1}{16}+\varepsilon}),
$$
where
$$
\Lambda_k(s)=(\frac{D^\ast|d|}{2\pi})^{s-1}\frac{\Gamma(s+k-1)}{\Gamma(k)}L_D(2s-1).
$$
So
$$
R_k =R_1
+\Lambda'_k(1)-\Lambda'_1(1)+O((D|d|)^{-\frac{1}{16}+\varepsilon}).
$$
As
$$
 \Lambda'_k(1)=\frac{\Gamma'(k)}{\Gamma(k)}L_D(1) +\log\frac{D^\ast|d|}{2\pi}L_D(1) + 2L'_D(1)\geq \Lambda'_1(1),
$$
we claim that
$$
R_1\geq  .0351,
$$
from which the estimate for $R_k$, which is stated at the
beginning of this section, follows. Write
$$
\frac{\zeta(s)L_D(s)}{\zeta(2s)} = \sum_{n=1}^{\infty}a_nn^{-s}
$$
with $a_1 = 1$, and $a_n \geq 0$. Then
$$
R_1= \sum_{n}a_{n}n^{-1} I(\frac{D^{\ast}|d|}{2\pi n^2}),
$$
where
$$
I(x) = \frac1{2\pi i} \int_{(2)} x^{s-1}
\frac{\Gamma(s)\zeta(4s-2)}{\zeta(2s-1)} \frac{ds}{(s-1)^2}.
$$
Miller-Yang [MY] proved that $I(x)>0$ if $x>0$ and that
$I(x)>.0351$ if $x\geq 4$. So we have $R_1\geq  .0351$ as claimed.

\enddocument